\documentclass[10 pt, leqno]{amsart}

\usepackage {amsfonts}
\usepackage{amsthm}
\usepackage{amssymb}
\usepackage{latexsym}
\usepackage{amsmath}

\theoremstyle{plain}
\newtheorem{tw}{Theorem}[section]

\newtheorem {lem} [tw]{Lemma}

\newtheorem{cor}[tw]{Corollary}

\theoremstyle{definition}
\newtheorem {deft}[tw] {Definition}
\newtheorem {rem} [tw]{Remark}

\newcommand{\bc} {\Bbb C}
\newcommand{\bn}{\Bbb N}

\newcommand{\bz}{\Bbb Z}

\newcommand{\alg} {\mathsf{A}}
\newcommand{\qg} {\mathsf{S}}
\newcommand{\QI} {\mathsf{\widetilde{QISO}^+}}
\newcommand{\qi}{\mathsf{QISO}^+}
\newcommand{\qil}{{\mathsf{QISO}}^{\mathcal L}}
\newcommand{\qib}{{\mathsf{QISO}}^{\bf Ban}}
\newcommand{\univ} {\mathsf{QU}}
\newcommand{\isox} {\mathsf{ISO}_X(X \times Y)}
\newcommand{\isoyx} {\mathsf{ISO}(X \times Y)}
\newcommand{\isoz} {\mathsf{ISO}(Z)}

\newcommand{\Alg} {\mathcal{A}}

\newcommand{\Hil}{\mathsf{H}}

\newcommand{\Ind}{\mathcal{J}}

\newenvironment{rlist}
{

\begin{enumerate}}
{\end{enumerate}}

\newcommand{\ot}{\otimes}
\newcommand{\ol}{\overline}

\newcommand{\wt}{\widetilde}

\numberwithin{equation}{section}


\begin{document}

\author{Jyotishman Bhowmick}
\footnote{\emph{Permanent address of the third named author:} Department of Mathematics, University of \L\'{o}d\'{z}, ul. Banacha 22, 90-238 \L\'{o}d\'{z}, Poland.}
\address{Stat-Math Unit, Indian Statistical Institute, 203, B. T. Road, Kolkata 700 208} \email{ jyotish\_r @isical.ac.in}
\author{Debashish Goswami}
\address{Stat-Math Unit, Indian Statistical Institute, 203, B. T. Road, Kolkata 700 208} \email{goswamid@isical.ac.in}
\author{Adam Skalski}
\address{Department of Mathematics and Statistics,  Lancaster University,
Lancaster, LA1 4YF} \email{a.skalski@lancaster.ac.uk}

\title{\bf Quantum isometry groups of $0$-dimensional manifolds}

\begin{abstract}
Quantum isometry groups of spectral triples associated with approximately finite-dimensional
$C^*$-algebras are shown to arise as inductive limits of quantum symmetry groups of corresponding
truncated Bratteli diagrams. This is used to determine explicitly the quantum isometry group of the
natural spectral triple on the algebra of continuous functions on the middle-third Cantor set. It
is also shown that the quantum symmetry groups of finite graphs or metric spaces coincide with the
quantum isometry groups of the corresponding classical objects equipped with natural Laplacians.
\end{abstract}

\maketitle


\section*{Introduction}
Following the pioneering ideas of Wang described in \cite{finsym} (and motivated by Connes), a number of mathematicians including
Bichon, Banica and others (\cite{graph}, \cite{metric}) have defined and studied the universal objects in certain categories of
quantum groups, all of which are quantum generalizations of groups acting on (typically finite) sets or algebras preserving some
given underlying structure like a metric or a functional. In this way, they have come up with several universal (compact) quantum
groups corresponding to the classical group of permutations (more generally, isometries w.r.t.\ a given metric) of a finite set,
the group of symmetries of a finite graph or the automorphism group of a finite dimensional matrix algebra. Motivated by their
work, the first two authors of the present article have begun a systematic effort to define and study similar universal quantum
groups beyond the finite-dimensional or `discrete' set-up, more precisely, in the framework of (possibly noncommutative)
differential geometry as proposed by Connes (\cite{connes}). They have been able to formulate a quantum group analogue of the
group of Riemannian isometries (\cite{Deb}, \cite{DebJyot}) as well as the group of orientation-preserving isometries
(\cite{Deb2}) of a (possibly noncommutative, given by spectral triple) Riemannian spin manifold. Many well-known and important
compact quantum groups (e.g. $SO_q(3)$) have been identified with such universal quantum groups for some suitably chosen spectral
triples.
            Recently in \cite{Chrivan} and  \cite{Chrivan2} Christensen and Ivan
constructed natural spectral triples on  approximately finite dimensional ($AF$) $C^*$-algebras.
The starting point for this article is the desire to understand and compute quantum isometry groups
of resulting noncommutative manifolds. $AF$ algebras  provide a  natural `connecting bridge'
between the finite and infinite dimensional noncommutative spaces and thus can be thought of as
$0$-dimensional manifolds. This is reinforced by the fact that Christensen and Ivan showed that on
each $AF$ algebra one can construct spectral triples with arbitrarily good summability properties.
We show that the `quantum group of orientation preserving isometries' of a Christensen-Ivan type
triple arises as an inductive limit of quantum isometry groups of certain finite-dimensional
triples (Theorem \ref{induc}). In the case when the $AF$ algebra in question is commutative, the
resulting quantum isometry groups of relevant finite-dimensional objects fit into the framework
described in the beginning of the introduction, as we show that they coincide with quantum symmetry
groups of finite graphs obtained by suitable truncations of the Bratteli diagrams. This observation
implies that the construction we consider can be thought of as giving a definition of a quantum
symmetry group of an arbitrary Bratteli diagram.  It also enables us to compute explicitly the
quantum isometry group of a spectral triple associated with the middle-third Cantor set introduced
first by Connes and later studied by Christensen and Ivan.  As a byproduct of our considerations of
this example we see that contrary to the classical case a quantum isometry of the product set
preserving the first factor in the suitable sense need not to be a product isometry.

    Having determined the universal objects for actions on $0$-dimensional noncommutative manifolds,
it is  natural to look back and see how one can accommodate the already existing theory of quantum
permutation and quantum automorphism groups of `finite' structures in the more general set-up of
quantum isometry groups. This is the second of the main objectives of the present article. We have
been able to identify the quantum group of automorphisms of a finite metric space or a finite graph
in the sense of Banica and Bichon with the quantum group of orientation (and suitable
`volume-form') preserving isometries of a natural spectral triple, thus successfully unifying the
approaches of \cite{metric} and \cite{graph} with that of \cite{Deb} and \cite{DebJyot}. We finish
the paper by suggesting a possible approach to defining quantum isometric actions on general
(compact) metric spaces and compute two explicit examples of universal quantum groups of
`isometries' in such context.

The detailed plan of the article is as follows: we begin by introducing basic notations and
recalling fundamental concepts related to quantum groups of (orientation preserving) isometries, as
defined in \cite{Deb} and \cite{Deb2}. Section 1 contains a description of the limit construction
for an inductive system of compact quantum groups and its application for quantum isometry groups.
In Section 2 we recall the construction of spectral triples on $AF$ algebras due to Christensen and
Ivan, show basic properties of the compact quantum groups appearing in the related inductive system
and relate them for commutative $AF$ algebras with the quantum symmetry groups of truncated
Bratteli diagrams. This is used to compute in Section 3 the quantum isometry group for Connes's
spectral triple related to the Cantor set. In Section 4 it is shown that the quantum symmetry group
of a finite metric space $X$ (\cite{metric}) coincides with the quantum isometry group resulting
from equipping the algebra of functions on $X$ with a natural Laplacian; the result has a natural
variant for the quantum symmetry group of a finite graph. Finally Section 5 contains a suggestion
of a tentative definition of  quantum isometry of a general metric space $X$ and computation of
such an object for $X=[0,1]$ and $X=S^1$.

\section*{Notations and preliminaries}

The symbol $\ot$ will always denote the minimal/spatial tensor product of $C^*$-algebras, purely algebraic tensor product will be denoted by $\odot$. We will occasionally use the language of Hilbert $C^*$-modules and multiplier algebras (see \cite{Lance}). Often when $X$ is a finite set we will write $C(X)$ to denote the algebra of all complex functions on $X$, with the point of view that when $X$ becomes an infinite topological space the correct generalisation is the algebra of continuous functions on $X$.

A \emph{compact quantum group} (c.q.g.) is a pair $(\qg, \Delta)$, where $\qg$ is a unital separable $C^*$-algebra and $\Delta : \qg \to \qg \ot \qg$ is a unital $C^*$-homomorphism satisfying the coassociativity: \\
(ai)  $(\Delta \ot {\rm id}) \circ \Delta=({\rm id} \ot \Delta) \circ \Delta$  \\
\hspace*{0.5 cm} and the quantum cancellation properties:\\
  (aii) the linear spans of $ \Delta(\qg)(\qg \ot 1)$ and $\Delta(\qg)(1 \ot \qg)$ are norm-dense in $\qg \ot \qg.$  \\
\hspace*{0.5 cm} Occasionally we will simply call $\qg$ a compact quantum group understanding by
this the existence of a suitable coproduct $\Delta$  on $\qg$. By a morphism in the category of
compact quantum groups we understand a unital $^*$-homomorphism intertwining the respective
coproducts.

A c.q.g.\ $(\qg,\Delta)$ is said  to \emph{(co)-act} on a unital $C^*$ algebra $\alg$ if there is
a unital $C^*$-homomorphism (called an action) $\alpha : \alg \to \alg \ot \qg$ satisfying the following conditions:\\
(bi) $(\alpha \ot {\rm id}) \circ \alpha=({\rm id } \ot \Delta) \circ \alpha$,\\
(bii) the linear span of $\alpha(\alg)(1 \ot \qg)$ is norm-dense in $\alg \ot \qg$.\\
\hspace*{0.5 cm} A \emph{unitary (co-) representation} of a compact quantum group $(\qg, \Delta)$
on a Hilbert space $\Hil$ is a linear map $U$ from $ \Hil $ to the $C^*$-Hilbert $\qg$-module $
\Hil \otimes \qg$ such that the  element $ \widetilde{U} \in {\mathcal M} ( {\mathcal K} ( \Hil )
\otimes \qg ) $ given by the formula $\widetilde{U}( \xi \ot b)=U(\xi)(1 \ot b)$ ($\xi \in \Hil, b
\in \qg)$) is a unitary satisfying
$$ ({\rm  id} \otimes \Delta ) \widetilde{U} = {\widetilde{U}}_{(12)} {\widetilde{U}}_{(13)}$$
In the last formula we used the standard `leg' notation: for an operator $X \in {\mathcal B}(\Hil_1
\ot \Hil_2),~  X_{(12)}$ and $X_{(13)}$ denote respectively the operators $X \ot I_{\Hil_2} \in
{\mathcal B}(\Hil_1 \ot \Hil_2 \ot \Hil_2)$ and $\sigma_{23} X_{12} \sigma_{23}\in {\mathcal
B}(\Hil_1 \ot \Hil_2 \ot \Hil_2)$ ($\sigma_{23}$ being the unitary on $\Hil_1 \ot \Hil_2 \ot
\Hil_2$ which flips the two copies of $\Hil_2$).

Given a unitary representation $U$ of $(\qg, \Delta)$ we  denote by $\alpha_U$ the
$^*$-homomorphism $\alpha_U(X)=\widetilde{U}(X \ot 1){\widetilde{U}}^*$ for $X \in {\mathcal
B}(\Hil)$. If $\tau$ is  a  not necessarily bounded, but densely  (in the weak operator topology)
defined linear functional on ${\mathcal B}(\Hil)$,  we say that $\alpha_U$ preserves $\tau$ if
$\alpha_U$ maps a suitable weakly dense $^*$-subalgebra   (say ${\mathcal D}$) in the domain of
$\tau$ into ${\mathcal D} \odot \qg$ and $( \tau \otimes {\rm id}) (\alpha_U(a))=\tau(a)1_\qg$  for
all $a \in {\mathcal D}$. When $\tau$ is bounded and normal, this is equivalent to the condition
$(\tau \ot {\rm id}) (\alpha_U(a))=\tau(a) 1_\qg$ being satisfied by all $a \in {\mathcal
B}(\Hil)$.

We say that a (possibly unbounded) operator $T$ on $\Hil$ commutes with $U$ if $T \ot I$ (with the
natural domain) commutes with $\widetilde{U}$. Sometimes such an operator will be called
\emph{$U$-equivariant}.

We briefly recall the definitions of quantum isometry groups, referring to \cite{Deb} and \cite{Deb2} for the details. Let
$({\mathcal A}^\infty, \Hil, D)$ be a spectral triple (of compact type, see \cite{Conspect}; note however that as we will often
consider here finite-dimensional objects we \emph{do not} require $D$ to be unbounded).  Consider the category ${\bf Q}^\prime$
whose objects are pairs $(\widetilde{\qg}, U)$, where $\widetilde{\qg}$ is a compact quantum group and $U$ is a unitary
representation of $\widetilde{\qg}$ in $\Hil$ such that the action $\alpha_U$ maps ${\mathcal A^{\infty}}$ into the ampliation of
its weak closure, and moreover $U$ commutes with $D$. The set of morphisms ${\rm Mor}((\qg,U),(\qg^\prime,U^\prime))$ is the set
of c.q.g.\ morphisms $\Phi : \qg \to \qg^\prime$ satisfying the condition $({\rm id} \ot \Phi) (U)=U^\prime$. If, additionally,
we are given a (possibly unbounded) positive operator $R$ such that $R$ commutes with $D$, we also consider the subcategory ${\bf
Q}^\prime_R$ of ${\bf Q}^\prime$ consisting of the $\qg$-actions for which $\alpha_U$ preserves the functional $\tau_R$ as in
\cite{Deb2}. It is proved in \cite{Deb2} that ${\bf Q}^\prime_R$ has a universal object, denoted by $\QI_R(\Alg^\infty, \Hil, D)$
or simply $\QI_R(D)$. We shall denote by $\qi_R(D)$ the  Woronowicz subalgebra  of $\QI_R(D)$ such that  $\alpha_U$ faithfully
maps  $ \Alg^\infty$ into $ (\Alg^\infty)^{''} \otimes \qi_R(D)$, and this subalgebra  is called the \emph{quantum group of
($R$-twisted) volume and orientation-preserving isometries} of the underlying spectral triple. It is also proved in \cite{Deb2},
Theorem 2.14  that under some further conditions a universal object in the bigger category ${\bf Q}^\prime$ exists. It is denoted
by $\QI(D)$. The corresponding  Woronowicz subalgebra for which $\alpha_U|_{\Alg^\infty}$ is faithful is denoted by $\qi(D)$, and
called the \emph{quantum group of orientation-preserving isometries}.

\section{Inductive limit construction for quantum isometry groups}

In this section we describe the limiting construction for an inductive system of compact quantum groups and give an application
for quantum isometry groups which is fundamental for the results of the next section.

The following lemma is probably known, but we include the proof for the sake of completeness.

\begin{lem}
\label{basic}
Suppose that $(\qg_n)_{n \in \bn}$ is a sequence of compact quantum groups and for each $n,m \in \bn,$ $n \leq m$ there is a c.q.g.\
morphism $\pi_{n,m}:\qg_n \to \qg_m$ with the compatibility property
\[ \pi_{m,k} \circ \pi_{n,m} = \pi_{n,k}, \;\;\; n\leq m \leq k.\]
Then the inductive limit of $C^*$-algebras $(\qg_n)_{n \in \bn}$ has a canonical structure of a compact quantum group. It will be
denoted $\qg_\infty$ or $\lim_{n \in \bn} \qg_n$. It has the following universality property:\\
for any c.q.g. $(\qg, \Delta)$ such that there are c.q.g. morphisms $\pi_n : \qg_n \to \qg$ satisfying for all $m,n \in \bn$, $m \geq n$ the equality $\pi_m \circ \pi_{n,m}=\pi_n$, there exists a unique c.q.g. morphism $\pi_\infty : \qg_\infty \to \qg$ such that $\pi_n=\pi_\infty \circ \pi_{n, \infty}$ for all $n\in \bn$, where we have denoted by $\pi_{n,\infty}$ the canonical unital $C^*$-homomorphism from $\qg_n$ into $\qg_\infty$.
\end{lem}

\begin{proof}
Let us denote the coproduct on $\qg_n$ by $\Delta_n$. We consider the unital $C^*$-homomorphism $\rho_n : \qg_n \to \qg_\infty
\ot \qg_\infty$ given by $ \rho_n=(\pi_{n,\infty} \ot \pi_{n,\infty}) \circ \Delta_n $, and observe that these maps  do satisfy
the compatibility property: $$ \rho_m \circ \pi_{n,m}=\rho_n~~~\forall n \leq m.$$ Thus, by the general properties of the
$C^*$-algebraic inductive limit, we have a unique unital $C^*$-homomorphism $\Delta_\infty : \qg_\infty \to \qg_\infty \ot
\qg_\infty$ satisfying $\Delta_\infty \circ \pi_{n,\infty}=\rho_n$ for all $n$. We claim that $(\qg_\infty, \Delta_\infty)$ is a
c.q.g.

We first  check that $ \Delta_{\infty} $ is coassociative.
It is enough to verify the coassociativity on the dense set $\cup_n \pi_{n,\infty}(\qg_n)$. Indeed, for $s=\pi_{n,\infty}(a)$ $(a \in \qg_n$),  by using $\Delta_\infty \circ \pi_{n,\infty}=(\pi_{n,\infty} \otimes \pi_{n,\infty}) \circ \Delta_n,$ we have the following:
\begin{align*}
(\Delta_\infty \otimes {\rm id})\Delta_\infty(\pi_{n,\infty}(a))\\
&= (\Delta_\infty \otimes {\rm id})(\pi_{n,\infty }\otimes \pi_{n,\infty})(\Delta_n(a))\\
&=(\pi_{n,\infty} \otimes \pi_{n,\infty} \otimes \pi_{n,\infty})(\Delta_n \otimes {\rm id})(\Delta_n(a))\\
&= (\pi_{n,\infty} \otimes \pi_{n,\infty} \otimes \pi_{n,\infty})( {\rm id} \otimes \Delta_n)(\Delta_n(a))\\
&=(\pi_{n,\infty} \otimes (\pi_{n,\infty} \otimes \pi_{n,\infty})\circ \Delta_n)(\Delta_n(a))\\
&=(\pi_{n,\infty} \otimes \Delta_\infty \circ \pi_{n,\infty})(\Delta_n(a))\\
&=({\rm id} \otimes \Delta_\infty)((\pi_{n,\infty} \otimes \pi_{n,\infty})(\Delta_n(a)))\\
&=({\rm id} \otimes \Delta_\infty)(\Delta_\infty(\pi_{n,\infty}(a))),
\end{align*} which proves the coassociativity.

Finally, we need to verify the quantum cancellation properties. Note that to show that  $ \Delta_{\infty} ( \qg_{\infty} ) ( 1 \otimes  \qg_{\infty} )  $ is dense in $ \qg_{\infty} \otimes  \qg_{\infty}$
it is enough to show that the above assertion is true with $\qg_{\infty}$ replaced by a dense subalgebra $ \bigcup_{n} \pi_{n, \infty} ( \qg_{n} ) $.

Using the density of $ \Delta_{n} ( \qg_{n} ) ( 1 \otimes \qg_{n} ) $ in $\qg_{n} \otimes \qg_{n} $ and the contractivity of the map $ \pi_{n, \infty} $ we note that $ ( \pi_{n, \infty} \otimes  \pi_{n, \infty} ) ( \Delta_{n} ( S_{n} ) ( 1 \otimes \qg_{n} ))$ is dense in $ ( \pi_{n, \infty} \otimes  \pi_{n, \infty} ) ( \qg_{n} \otimes \qg_{n} ).$
 This implies that  $ ( \pi_{n, \infty} \otimes  \pi_{n, \infty} ) ( \Delta_{n} ( \qg_{n} ) ) ( 1 \otimes   \pi_{n, \infty} ( \qg_{n} ) ) $ is dense in $  \pi_{n, \infty} ( \qg_{n} ) \otimes \pi_{n, \infty} ( \qg_{n} ) $ and hence $ \Delta_{\infty} ( \pi_{n, \infty} ( \qg_{n} ) ) ( 1 \otimes \pi_{n, \infty} ( \qg_{n} ) ) $ is dense in  $ \pi_{n, \infty} ( \qg_{n} ) \otimes  \pi_{n, \infty} ( \qg_{n} ).$ The proof of the claim now follows by  noting that $ \pi_{n, \infty} ( \qg_{n} ) =  \pi_{m, \infty} \pi_{n, m} ( \qg_{n} ) \subseteq \pi_{m, \infty} ( \qg_{m} ) $ for any $m \geq n$, along with the above observations. The right quantum cancellation property can be shown in the same way.

The proof of the universality property is routine and hence omitted.

\end{proof}

Note that the proof remains valid for any other indexing set for the net, not necessarily $\bn$.

The next theorem connects the inductive construction above with some specific quantum isometry groups.

\begin{tw} \label{induc}
Suppose that $\alg$ is a $C^*$-algebra acting on a Hilbert space $\Hil$ and that  $D$ is a (densely defined) selfadjoint operator
on $\Hil$ with compact resolvent, such that $D$ has a one-dimensional eigenspace spanned by a vector $\xi$ which is cyclic and
separating for $\alg$. Let $(\Alg_n)_{n \in \bn}$ be an increasing net of a unital $^*$-subalgebras of $\alg$ and put $\Alg =
\bigcup_{n \in \bn} \Alg_n$. Suppose that $\Alg$ is dense in $\alg$ and that for each $a \in \Alg$ the commutator $[D,a]$ is
densely defined and bounded. Additionally put $\Hil_n = \ol{\Alg_n \xi}$, let $P_n$ denote the orthogonal projection on $\Hil_n$
and assume that each $P_n$ commutes with $D$. Then each $(\Alg_n, \Hil_n, D|_{\Hil_n})$ is a spectral triple satisfying the
conditions of Theorem 2.14 of \cite{Deb2}, there exist natural compatible   c.q.g.\ morphisms $\pi_{m,n}:\QI(\Alg_m, \Hil_m,
D|_{\Hil_m}) \to \QI(\Alg_n, \Hil_n, D|_{\Hil_n})$ ($n, m \in \bn, m \leq n)$ and
\[ \QI(\Alg, \Hil, D) = \lim_{n \in \bn} \QI(\Alg_n, \Hil_n, D|_{\Hil_n}). \]
Similar conclusions hold if we replace everywhere above $\QI$ by $\qi$.
\end{tw}
\begin{proof}
We prove the assertion corresponding to $\QI$ only, since the proof for $\qi$  follows by very similar arguments.
 Let us denote $ \QI ( \Alg_{n}, \Hil_{n}, D_{n} ) $ by $ \qg_{n} $ and the corresponding unitary representation (in $\Hil_n$) by  $ U_{n}.$
Let us denote the category of compact quantum groups acting by orientation preserving isometries on
$ (\Alg_{n}, \Hil_n, D|_{\Hil_n}) $ and $ (\Alg, \Hil, D) $ respectively by $ \bf C_{n} $ and $ \bf
C $ .

Since $ U_{n} $ is a unitary which commutes with $ D_n \equiv D|_{\Hil_n} $ and hence preserves the
eigenspaces of $ D_n$,  it restricts to a unitary representation of $ S_{n} $ on each $ H_{m} $ for
$ m \leq n.$ In other words, $(\qg_n, U_n|_{\Hil_m}) \in {\rm Obj}( {\bf C_m} )$, and by
 the universality of $ \qg_{m}$ there exists a compact quantum group morphism, say, $ \pi_{m,n} : \qg_{m} \rightarrow \qg_{n} $ such that $ ( {\rm id} \otimes \pi_{m,n} ) U_{m}|_{\Hil_{m}} = U_{n}|_{\Hil_{m}}.$

Let $ p \le m \le n .$ Then  we have $ ( {\rm id}  \otimes \pi_{m,n} \pi_{p,m} )U_{p}|_{\Hil_{p}} =
U_{n}|_{\Hil_{p}} .$ It follows  by the uniqueness of the map $ \pi_{p,n}$ that $ \pi_{p,n} =
 \pi_{m,n}\pi_{p,m} $, i.e.\   $(\qg_n)_{n \in \bn}$ forms an inductive system of compact quantum
groups satisfying the assumptions of  Lemma \ref{basic}.  Denote by $\qg_\infty$ the inductive
limit c.q.g.\ obtained in that lemma, with $\pi_{n,\infty} : \qg_n \to \qg$ denoting the
corresponding c.q.g.\ morphisms. The family of formulas $U|_{\Hil_n}:=({\rm id} \ot \pi_{n,\infty})
\circ U_n$ combine to define a unitary representation $U$ of $\qg_\infty$ on $\Hil$.  It is also
easy to see from the construction that $U$ commutes with $D$. This means that $(\qg_\infty, U) \in
{\rm Obj}({\bf C})$, hence there exists a unique surjective c.q.g.\ morphism from $\qg:=\QI(\Alg,
\Hil, D)$ to $\qg_\infty$ identifying $\qg_\infty$ as a quantum subgroup of $\qg$.

The proof will now be complete if we can show that there is a surjective c.q.g.\ morphism in the
reverse direction, identifying $\qg$ as a quantum subgroup of $\qg_\infty$. This can be deduced
from Lemma \ref{basic} by using the universality property of the inductive limit. Indeed,
 for each $n \in \bn$ the unitary representation, say $V_n$, of $\QI(\Alg, \Hil, D)$ restricts to $\Hil_n$ and commutes with $D$ on that subspace, thus inducing a c.q.g.\ morphism $\rho_n$  from $\qg_n=\QI(\Alg_n,\Hil_n, D_n)$ into $\qg$. The family of morphisms $(\rho_n)_{n \in \bn}$ satisfies the compatibility conditions required in Lemma \ref{basic}. It remains to show  that the induced c.q.g.\ morphism $\rho_{\infty}$  from $\qg_\infty$ into $\qg$ is surjective.
By the faithfulness of the representation $V$ of $\QI(\Alg, \Hil, D)$, we know that the span of matrix elements corresponding to all $V_n$  forms a norm-dense subset of $\qg$.  As the range of $\rho_n$ contains the matrix elements corresponding to $V_n=V|_{\Hil_n}$, the proof of surjectivity of $\rho_{\infty}$ is finished.
\end{proof}

The assumptions of the theorem might seem very restrictive. In the next section however we will describe a natural family of spectral triples on $AF$-algebras, constructed in \cite{Chrivan}, for which we have exactly the situation as above.

\section{Quantum isometry groups for spectral triples on AF algebras}

We first recall the construction of natural spectral triples on $AF$ algebras due to
E.\,Christensen and C.\,Ivan (\cite{Chrivan}). Let $\alg$ be a unital $AF$ $C^*$-algebra, the norm
closure of an increasing sequence $(\alg_n)_{n \in \bn}$ of finite dimensional $C^*$-algebras. We
always put $\alg_0 = \bc 1_{\alg}$, $\Alg = \bigcup_{n=1}^{\infty} \alg_n$  and assume that the
unit in each $\alg_n$ is the unit of $\alg$.
 Suppose that $\alg$ is acting on a Hilbert space $\Hil$ and that $\xi \in \Hil$ is a separating and cyclic unit vector for $\alg$.
Let $P_n$ denote the orthogonal projection onto the subspace $\Hil_n:=\alg_n \xi$ of $\Hil$ and
write $Q_0=P_0=P_{\bc \xi}$, $Q_n=P_n - P_{n-1}$ for $n \in \bn$. There exists a (strictly
increasing) sequence of real numbers $(\alpha_n)_{n=1}^{\infty}$ such that the selfadjoint operator
$D=\sum_{n \in \bn} \alpha_n Q_n$ yields a spectral triple $(\Alg, \Hil,D)$ such that the topology
on the state space of $\alg$ induced by the Rieffel metric (\cite{Riefstate}) coincides with the
weak$^*$-topology. Due to the existence of a cyclic and separating vector the orientation
preserving quantum isometry group exists by Theorem 2.14 of \cite{Deb2}.

In \cite{Chrivan} it was additionally observed that if $\alg$ is infinite-dimensional and $p>0$
then one can choose $(\alpha_n)_{n=1}^{\infty}$ in such a way that the resulting Fredholm module is
$p$-summable (\cite{Conspect}). This reflects the fact that $AF$ algebras should be thought of as
$0$-dimensional noncommutative spaces.

Note that for each $n\in\bn$ by restricting we obtain a (finite-dimensional) spectral triple
$(\alg_n, \Hil_n, D|_{\Hil_n})$. As we are precisely in the framework of Theorem \ref{induc}, to
compute $\qi(\Alg,\Hil,D)$  we need to understand the quantum isometry groups $\qi(\Alg_n, \Hil_n,
D|_{\Hil_n})$ and embeddings relating them. To simplify the notation we will write
$\qg_n:=\qi(\alg_n, \Hil_n, D|_{\Hil_n})$.

We begin with some general observations.

\begin{lem}
Let  $\univ_{\alg_n, \omega_{\xi}}$ denote the universal quantum group acting on $\alg_n$ and
preserving the (faithful) state on $\alg_n$ given by vector $\xi$ (see \cite{finsym}). There exists
a  c.q.g.\ morphism from $\univ_{\alg_n, \omega_{\xi}}$ to $\qg_n$.
\end{lem}

\begin{proof}
The proof is based on considering the spectral triple given by $(\alg_n, \Hil_n, D'_n)$, where $D'_n = P_n - P_0$. It is then
easy to see that $\qi(\alg_n, \Hil_n, D'_n)$ is isomorphic to the universal compact quantum group acting on $\alg_n$ and
preserving $\omega_{\xi}$. On the other hand universality assures the existence of the   c.q.g.\ morphism from
$\qi(\alg_n, \Hil_n, D'_n)$ to $\qg_n$.
\end{proof}


\begin{lem} \label{univ}
Assume that each $\alg_n$ is commutative, $\alg_n = \bc^{k_n}$, $n \in \bn$. There exists a
c.q.g.\ morphism from $\univ_{k_n}$ to $\qg_{k_n}$, where $\univ_{k_n}$ denotes  the universal
quantum group acting on $k_n$ points (\cite{finsym}).
\end{lem}

\begin{proof}
 The proof is identical to the one above - we only need to observe additionally that for any measure $\mu$ on the set $\{1,\ldots,k_n\}$ which has   full support there is a natural c.q.g.\ morphism from $\univ_{k_n}$ to $\univ_{\bc^{k_n}, \mu}$. In case when $\mu$ is uniformly distributed, we simply have $\univ_{\bc^{k_n}, \mu} = \univ_{k_n} $, as follows from the first part  of Lemma \ref{fact} below.
\end{proof}

Let $\alpha_n: \alg_n \to \alg_n \ot \qg_n$ denote the universal action (on the $n$-th level). Then we have the following
important property, being the direct consequence of the Theorem \ref{induc}. We have
\begin{equation} \label{invar}\alpha_{n+1} (\alg_n)
\subset \alg_n \ot \qg_{n+1}\end{equation} (where we identified $\alg_n$ with a subalgebra of $\alg_{n+1}$) and $\qg_n$ is
generated exactly by these coefficients of $\qg_{n+1}$ which appear in the image of $\alg_n$ under $\alpha_{n+1}$. This in
conjunction with the previous lemma suggests the strategy for computing relevant quantum isometry groups inductively. Suppose
that we have determined the generators of $\qg_n$. Then $\qg_{n+1}$ is generated by generators of $\qg_n$ and these of the
$\univ_{\alg_n, \omega_{\xi}}$, with the only additional relations provided by the equation \eqref{invar}.

This will be used below to determine the concrete form of relations determining $\qg_n$ for the commutative AF algebras.

\begin{lem} \label{embed}
Let $\alg$ be a commutative AF algebra. Suppose that $\alg_n$ is isomorphic to $\bc^{m}$ and the embedding of $\alg_{n}$ into
$\alg_{n+1}$ is given by a sequence $(l_i)_{i=1}^{m}$. Let $m'= \sum_{i=1}^{m} l_i$. Suppose that the `copy' of $\univ_{m}$ in
$\qg_n$ is given by the family of projections $a_{i,j}$ ( $i,j\in\{1, \ldots m\}$ ) and that the `copy' of $\univ_{m'}$ in
$\qg_{n+1}$ is given by the family of projections $a_{(i,r_i),(j,s_j)}$ ($i,j\in \{1,\ldots,m\},$ $r_i\in\{1, \ldots, l_i\} $,
$s_j \in \{1, \ldots,l_j\}$). Then the formula \eqref{invar} is equivalent to the following system of equalities:
\begin{equation} \label{finite} a_{i,j} = \sum_{r_{i}=1}^{l_i} a_{(i,r_i), (j,s_j)}\end{equation}
for each $i,j \in \{1, \ldots,m\},s_{j}\in \{1,\ldots, l_j\}$.
\end{lem}
\begin{proof}
We have (for the universal action $\alpha: \alg_n \to \alg_n \ot \qg_n$)
\[\alpha(\wt{e_i} ) = \sum_{j=1}^m \wt{e_j} \ot a_{i,j}\]
where by $\wt{e_i}$ we denote the image of the basis vector $e_i\in \alg_n$ in $\alg_{n+1}$. As $\wt{e_j} = \sum_{r_j=1}^{l_j}
e_{(j,s_j)}$,
\[\alpha(\wt{e_i}) = \sum_{r_i=1}^{l_i} \alpha(e_{i,r_i}) = \sum_{r_i=1}^{l_i} \sum_{j=1}^m \sum_{s_j=1}^{l_j}
e_{(j, s_j)}  \ot a_{(i,r_i), (j,s_j)}.\] On the other hand we have
\[\alpha(\wt{e_i} ) = \sum_{j=1}^m \sum_{s_j=1}^{l_j}
e_{(j,s_j)}  \ot a_{i,j},\] and the comparison of the formulas above yields exactly \eqref{finite}.
\end{proof}

One can deduce from the above lemma the exact structure of generators and relations between them for each $\qg_n$ associated with
a commutative AF algebra. To be precise, if $\alg_n= \bc^{k_n}$ for some $k_n \in \bn$, then the quantum isometry group $\qg_n$
is generated as a unital $C^*$-algebra by the family of selfadjoint projections $\bigcup_{i=1}^{n} \{a_{\alpha_i,
\beta_i}:\alpha_i, \beta_i = 1,\cdots, k_i\}$ such that for each fixed $i=1,\ldots,n$ the family $ \{a_{(\alpha_i,
\beta_i)}:\alpha_i, \beta_i = 1,\cdots, k_i\}$ satisfies the relations of $\univ_{k_n}$ and the additional relations between
$a_{(\alpha_i, \beta_i)}$ and $a_{(\alpha_{i+1}, \beta_{i+1})}$ for $i \in \{1,\ldots,n-1\}$ are given by the formulas
\eqref{finite}, after suitable reinterpretation of indices according to the multiplicities in the embedding of $\bc^{k_i}$ into
$\bc^{k_{i+1}}$.

In \cite{graph} J.\,Bichon introduced the notion of a quantum symmetry group of a finite directed
graph. As each AF algebra can be described via its Bratteli diagram, it is natural to ask whether
the construction in this paper can be compared to the one in \cite{graph}. We begin by stating some
elementary facts in the following lemma.

\begin{lem}
\label{fact}
 Let $ \alpha $ be an action of a c.q.g.\ $ \qg $ on $C(X)$ where $ X $ is a finite set.  Then $\alpha$ automatically preserves the functional $\tau$   corresponding to the counting measure:
$$ ( \tau \otimes {\rm id} )( \alpha ( f )) = \tau ( f ).1_{\qg}.$$  Thus $ \alpha $ induces a unitary $ \widetilde{\alpha} \in {\mathcal B }( l^{2} ( X ) ) \otimes \qg $ given by $ \widetilde{\alpha} ( f \otimes q ) = \alpha ( f ) ( 1 \otimes q ).$
If we define $ \alpha^{(2)} : C ( X ) \otimes C (X) \rightarrow C (X ) \otimes C ( X ) \otimes \qg$ by $\alpha^{(2)} = ( {\rm id}_{2} \otimes m_\qg ) \sigma_{23} ( \alpha \otimes \alpha ),$ where $m_\qg$ denotes the multiplication map from $\qg \otimes \qg$ to $\qg$, and ${\rm id}_2$ denotes the identity map on $C(X) \otimes C(X)$, then $ \alpha^{(2)} $ leaves the diagonal algebra $ C ( D_{ X \times X } ) $ invariant (here  $ D_{X \times X} = \{ ( x,x ) : x \in X \}$).
   \end{lem}
\begin{proof}
    Let $X=\{ 1,...,n\}$ for some $n \in \bn$ and denote by $\delta_i$ the characteristic function of the point $i$. Let $\alpha(\delta_i)=\sum_j \delta_j \ot q_{ij}$ where $ \{q_{ij}:i,j=1\ldots n\}$  are the images of the canonical  generators of the quantum permutation group as in \cite{finsym}. Then  $\tau$-preservation of $\alpha$ follows from the properties of the generators of the quantum permutation group, which in particular imply that $\sum_j q_{ij}=1=\sum_i q_{ij}$.

Using the fact that $q_{ij}$ and $q_{ik}$ are orthogonal for $i,j,k \in \{1,\ldots,n\}$, $ j \neq k$, we obtain    \[ \alpha^{(2)} ( \delta_{i} \otimes \delta_{i} )
=  \sum_{j,k} \delta_{j} \otimes \delta_{k} \otimes q_{ij} q_{ik}
=  \sum_{j} \delta_{j} \otimes \delta_{j} \otimes q_{ij},\] from which the invariance of the diagonal under $\alpha^{(2)}$ is immediate.

The other statements follow easily.
        \end{proof}

 Observe that  since $ \widetilde{\alpha^{(2)}} $ is a unitary, it leaves the $\tau$-orthogonal complement of $ C ( D_{X \times X} ) $ in $l^2(X \times X)$, i.e,\  the space of functions on the set $Y=\{ ( x,y ) : x,y \in X, x \neq y \}  $ invariant as well.

Recall now the definition of the quantum automorphism group of a finite graph given in \cite{graph}. Let $ ( V, E ) $ be a graph
with $V$ denoting the set of vertices and $E$ the set of edges. Let $ s : E \rightarrow V $ ( respectively $ t : E \rightarrow V
$ ) be the source map ( respectively the target map ). The target and source maps induce $^*$-homomorphisms $ s^{*}, t^{*} : C (
V ) \rightarrow C ( E ).$ Let $ m : C ( E ) \otimes C ( E ) \rightarrow C ( E ) $ be the pointwise multiplication map on $E$ and
given a quantum group action $\alpha$ on $C(V)$ let $ \alpha^{(2)} $ be defined  as in Lemma \ref{fact}.

\begin{deft}{(\cite{graph})}
An action of a c.q.g.\ $ \qg $ on a finite graph $ G = ( V, E ) $ consists of  an action $\alpha$
of $ \qg $ on the set of vertices, $ \alpha : C ( V ) \rightarrow C ( V ) \otimes \qg $ and an
action $\beta$ of $ \qg $ on the set of edges, $ \beta : C ( E ) \rightarrow  C ( E ) \otimes \qg
$, such that \[ ( ( m ( s_{*} \otimes t_{*} ) ) \otimes \textup{id}_{\qg}) \circ \alpha^{(2)} =  \beta\circ  ( m (
s^{*} \otimes t^{*} ) ).\] It is also called a quantum symmetry of the graph $(V,E)$.
\end{deft}

The \emph{quantum automorphism group of the finite graph} is the universal object in the category
of compact quantum groups with actions as above. We refer to \cite{graph} for the details.

Let us now restrict our attention to a truncation of a Bratteli diagram (up to $n$-the level, say),
of a commutative $AF$ algebra. The set of vertices $V$ is a disjoint union of sets $V_1,..., V_n$
with $V_1$ being a singleton, and there exist surjective maps $\pi_j : V_j \to V_{j-1}$ ($j \geq
2$) determining the graph structure. Denote by $\pi$ the map from $V$ to $V$ defined by the
formulas $\pi|_{V_j}=\pi_j$ for $j \geq 2$, and $\pi={\rm id}$ on $V_1$. Then $\pi^* : C(V) \to
C(V)$ is a $C^*$-homomorphism, with $\pi^*|_{C(V_j)}$ injective for each $j \leq n-1$. The
corresponding graph is obtained by joining $\pi_{i+1}(v)\in V_i$ with $v\in V_{i+1}$ for each
$i=1,\ldots, n$, $v \in V_i$.

Denote by $ D ( V ) $ the diagonal subalgebra of $ C ( V ) \otimes C ( V ),$ i.e.\ the span of $ \{
\delta_{v} \otimes \delta_{v} : v \in V \}.$ Since the map $ m (  s_{*} \otimes t_{*} ) $ is onto,
$ C ( E ) \cong ( C ( V ) \otimes C ( V ) )/ \textup{Ker} ( m ( s_{*} \otimes t_{*} ) ).$ Indeed,
for the graph corresponding to the commutative $AF$ algebra described above,  $ C ( E ) $ is
isomorphic to the subalgebra $ ({\rm  id} \otimes \pi^{*} ) (D (V))  $ of $ C ( V ) \otimes C ( V )
,$ and the condition on $\beta$ in the above definition of a quantum symmetry of the graph $(V,E)$
amounts to saying that $\alpha^{(2)}$ leaves the algebra ${\mathcal C}:=({\rm  id} \otimes \pi^{*}
) (D (V))$ invariant, in which case $\beta$ is, up to the obvious identification, nothing but the
restriction of $\alpha^{(2)}$ on ${\mathcal C}$.

In other words, an equivalent description of the objects in the category of the quantum symmetry of
such a finite commutative (i.e.\ with all matrix algebras in the vertices being one-dimensional)
Bratteli diagram is obtained via observing that they correspond precisely to these  c.q.g.\ actions
$\alpha$ on $C(V)$ for which ${\mathcal C}$ is left invariant by $\alpha^{(2)}$. This leads to the
following result:

\begin{tw} \label{AFgraph}
Let $\alg$ be a commutative $AF$ algebra. Then the quantum isometry group $\qg_n$ described in the beginning of this section coincides with the quantum symmetry group of the graph given by the restriction of the Bratteli diagram of $\alg$ to the $n$-th level.
\end{tw}

\begin{proof}

Suppose first that we are given a quantum isometry of the canonical spectral triple on the
respective `finite' part of the $AF$ algebra in question, so that we have a c.q.g.\ action $(\qg,
\alpha)$ on $C(V)$ such that each $\alpha_j \equiv \alpha|_{C(V_j)}$ leaves  $C(V_j)$ invariant
($j=1,\ldots,n$) and that $\alpha$ commutes with the embeddings $\pi_j$, that is
\[ \alpha_{j + 1} {\pi_{j+1}}^{*} = ( {\pi_{j+1}}^{*} \otimes \textup{id} ) \alpha_{j }.\]
We deduce that
\begin{align*}{\alpha}^{(2)}_{j + 1} ( \textup{id} \otimes {\pi_{j+1}}^{*} ) &= ( \textup{id}_{2} \otimes m_\qg ) \sigma_{23} ( \alpha_{j + 1} \otimes \alpha_{j + 1} {\pi_{j+1}}^{*} ) \\& = ( id_{2} \otimes m_\qg ) \sigma_{23} ( \alpha_{j + 1} \otimes ( {\pi_{j+1}}^{*} \otimes {\rm id} ) \alpha_{j }  ) = ( {\pi_{j+1}}^{*} \otimes{\rm  id} ) {\alpha}^{(2)}_{j} \end{align*}
 Using the above expression and the fact that $ {\alpha}^{(2)} $ leaves $D(V)$ invariant (by the second part of Lemma \ref{fact}), we see that $ \alpha^{(2)} $ leaves ${\mathcal C}$ invariant.

   Conversely, we need to show that a quantum action of a c.q.g.\ on the Bratteli diagram
induces a quantum symmetry of  the corresponding part of a spectral triple on the $AF$ algebra. Let
$(\qg, \alpha)$ be then an action on $C(V)$ such that $\alpha^{(2)}$ leaves ${\mathcal C}$
invariant. It follows from the discussion before the lemma that we have the corresponding action
$\beta$ on $C(E)$. Therefore we can start with an action $ \alpha $ on the Bratteli diagram such
that $ {\alpha^{(2)}}_{C(V_n)} $ preserves $ ( {\rm id} \otimes {\pi_{n}}^{*} ) ( D (V_{n -1}) ).$  We first
show by induction that $\alpha$ leaves each $C(V_j)$ invariant. Consider $j=n$ first. Observe that
$C(V_n)$ is nothing but  $ \textup{Ker}(\Psi)$, where $\Psi : C(V) \to C(E)$ is the map $f \mapsto
( m \circ ( s^{*} \otimes t^{*} ))(f \ot 1)$. It is clear from the definition of a quantum symmetry
of a graph that $\alpha$ will leave this subalgebra invariant. This implies that $\alpha$ (being a
unitary w.r.t.\ the counting measure on $C(V)$) will leave $C(V_1 \cup...\cup V_{n-1})$ invariant
as well, and thus restricts to a quantum symmetry of the reduced graph obtained by deleting $V_n$
and the corresponding edges. Then the inductive arguments complete the proof that each $C(V_j)$ is
left invariant.

The proof will now be complete if we can show  that
\begin{equation} \label{alcom} \alpha_{m+ 1}
{\pi^*_{m+1}} = ( {\pi^*_{m+1}}\otimes{\rm  id} )
 \alpha_{m } \end{equation}
 for each $m=1,\ldots,n-1$. Let $V_m=\{ v^m_1, \ldots,v^m_{t_m} \},$ and let $q_{m,lj}$ be elements of $\qg$ such that $$ \alpha ( \delta_{ v^{m}_{l}} ) = \sum_{j}  \delta_{ v^{m}_{j}} \otimes q_{m, lj}. $$
We set  $ \Lambda_j = \{ m : \pi( m ) = j \} .$
Then we have \begin{align*}
 \alpha^{(2)} ( \textup{id} \otimes  {\pi^*_{m+1}} ) ( \delta_{v^{m}_{i}} \otimes  \delta_{ v^{m}_{i}} )
 &= \alpha^{(2)} (  \delta_{ v^{m}_{i}} \otimes \sum_{k \in \Lambda_{i}}  \delta_{ v^{m + 1}_{k}} )
 \\&= \sum_{j,p}  \sum_{k \in \Lambda_{i} }  \delta_{ v^{m}_{j}} \otimes  \delta_{ v^{m + 1}_{p}} \otimes q_{m, ij} q_{m+ 1, kp}. \end{align*}
Further then for all  $p \notin \Lambda_{j}$
\[ \sum_{k \in \Lambda_{i}  } q_{m, ij} q_{m + 1, kp} = 0. \]
   Multiplying by $ q_{m + 1, ~ k^{'}p}$ where $ k^{'} \in \Lambda_{i},$ we obtain that
\begin{equation} \forall_{p \notin \Lambda_{j}, k \in \Lambda_{i}}\;\; q_{m, ij} q_{m + 1, kp} = 0 \label{gr1}\end{equation}
As stated after Lemma \ref{fact}, $\alpha^{(2)}$  leaves the ortho-complement of the diagonal algebra in $l^2(V \times V)$ invariant. This means that if $ k \notin \Lambda_{i}$ then $ \alpha^{(2)} ( \delta_{v^{m}_{i}} \otimes \delta_{ v^{m+ 1}_{k}} )$ belongs to $  C ( D^{c} ) \otimes \qg $. On the other hand
$$ \alpha^{(2)} ( \delta_{ v^{m}_{i}} \otimes \delta_{ v^{m + 1 }_{k}} ) = \sum_{j,p} \delta_{ v^{m}_{j}} \otimes \delta_{ v^{m + 1 }_{p}} \otimes q_{m, ij} q_{m + 1, kp},$$
so  \begin{equation} \forall_{k \notin \Lambda_{i}, p \in  \Lambda_{j}} \;\;\; q_{m, ij} q_{m + 1, kp} = 0.
\label{gr2}\end{equation} Further
$$ \alpha_{m+1} {\pi^*_{m+1}} ( \delta_{ v^{m}_{i} }) =
\sum_{j} \delta_{v^{m + 1}_{j}} \otimes \sum_{k \in \Lambda_{i}} q_{m + 1, kj}.$$ We also have  $$
( {\pi^*_{m+1}} \otimes \textup{id} ) \alpha_{m} ( \delta_{ v^{m}_{i} }) = \sum_{j} \delta_{v^{m +
1}_{j}} \otimes \sum_{r:~ j \in \Lambda_{r}} q_{m , ir} = \sum_{j}  \delta_{v^{m + 1}_{j}} \otimes
q_{m, i \pi(j)}.$$
Finally  \begin{align*} \sum_{k \in \Lambda_{i}} q_{m + 1, kp} &= \sum_{k \in
\Lambda_{i} } ( \sum_{j} q_{m, ij} )q_{m+ 1, kp}
 = \sum_{k \in \Lambda_{i} } ( \sum_{j} q_{m, ij} q_{m+ 1, kp} ) \\&= \sum_{k \in \Lambda_{i} } q_{m, i \pi ( p )} q_{m + 1,kp} = \sum_{k} q_{m, i \pi ( p )} q_{m+ 1,kp} - \sum_{k \notin \Lambda_{i}} q_{m, i \pi ( p )} q_{m + 1,kp} \\&=  q_{m, i \pi ( p )} ,\end{align*}
where in the third equality we used \eqref{gr1} and in the final equality we used \eqref{gr2} as well as the relation $\sum_k q_{m+1,kp}=1$. This
 shows that \eqref{alcom} holds 
and the proof is finished.

\end{proof}
%

The above result justifies the statement that the quantum isometry groups of Christensen-Ivan triples on $AF$ algebras provide
natural notions of quantum symmetry groups of the corresponding Bratteli diagrams. The theorem could be proved directly by
comparing  the commutation relations of \cite{graph} with these listed in Lemma \ref{embed}, but the method we gave has the
advantage of being more functorial and transparent.

\section{The example for the middle-third Cantor set}

In the special case when $\alg$ is the  (commutative) AF algebra of continuous functions on the
middle-third Cantor set we can use the observations of the last section to provide explicit
description of the quantum isometry groups $\qg_n$, and therefore also of $\QI(\Alg, \Hil, D)$.
Note that several variants of the  spectral triple we consider here were studied in \cite{Chrivan}
and in \cite{Chrivan2}, where its construction is attributed to the unpublished work of Connes.

\begin{tw} \label{mid3}
Let $\alg$ be the $AF$ algebra arising as a limit of the unital embeddings
\[\bc^2 \longrightarrow \bc^2 \ot \bc^2 \longrightarrow \bc^2 \ot \bc^2 \ot \bc^2 \longrightarrow \,\cdots.\]
Suppose that the state $\omega_{\xi}$ is the canonical trace on $\alg$.
Then $\qg_1 = C(\bz_2)$ and for $n \in \bn$
\[ \qg_{n+1} = (\qg_n \star \qg_n) \oplus (\qg_n \star \qg_n).\]
\end{tw}
\begin{proof}
Begin by noticing that for each $n \in \bn$ we have $\alg_n =\bc^{2^n}$ and therefore  we can use a natural multi-index notation for the indexing sets discussed in the paragraph after Lemma \ref{embed}. Precisely speaking denote for each $n \in \bn$ by $\Ind_n$ the set
$\{i_1i_2\cdots i_n: i_j\in\{1,2\} \textrm{ for }j=1,\ldots,n\}$ . Multi-indices in $\Ind:= \bigcup_{n \in \bn} \Ind_n$ will be denoted by capital letters $I,J,\ldots$ and let the basis of the algebra $\alg_n$ be indexed by elements of $\Ind_n$. Then the natural embeddings between
$\alg_n$ and $\alg_{n+1}$ can be conveniently described by the formula
\[e_I \longrightarrow e_{I1} + e_{I2},\;\;\; I \in \Ind_n,\]
where we use the standard concatenation of multi-indices. The equations \eqref{finite} take now the following form:
\begin{equation} \label{2finite}  a_{I,J} = a_{I1,J1} +a_{I2,J1} = a_{I1,J2} + a_{I2,J2}, \;\; I,J \in \Ind.\end{equation}
 Note that so far the fact that at every step we `divide the set' into $2$ parts did not play any significant role, we could adopt the multi-index notation replacing $2$ by $3,4,\ldots$. Analyse now the equations \eqref{2finite} remembering that both $\{a_{I,J}:I,J \in \Ind_n\}$ and
$\{a_{I',J'}:I',J' \in \Ind_{n+1}\}$ form `magic unitaries' whose entries are orthogonal
projections (\cite{finsym}). Fix $I \in \Ind_n$ and consider the respective equalities
\[1 = \sum_{J\in \Ind_n} a_{I,J} = \sum_{J \in \Ind_n} (a_{I1,J1} + a_{I2,J1})\]
and
\[ 1 = \sum_{J' \in \Ind_{n+1}} a_{I1, J'} = \sum_{J \in \Ind_{n}} (a_{I1, J1} + a_{I1,J2}) .\]
They imply that
\[ \sum_{J \in \Ind_n} a_{I2,J1} =  \sum_{J \in \Ind_{n}}  a_{I1,J2} .\]
By formulas \eqref{2finite} each of the respective factors in the sum above is an orthogonal
subprojection of the projection $a_{I,J}$ and the projections $a_{I,J}$, $a_{I,K}$ are mutually
orthogonal when $J\neq K$. Therefore  we actually must have
\[  a_{I2,J1} =  a_{I1,J2} \]
for each $J \in \Ind_n$. Using once again \eqref{2finite} we see that actually
\begin{equation} \label{2final} a_{I1,J1}=a_{I2,J2} = a_{I,J} - a_{I1,J2} = a_{I,J} - a_{I2,J1}, \;\;\ I,J \in \Ind_n. \end{equation}
This means that the choice of $a_{I1,J1}$ determines already the remaining three projections. Observe also that if we fix the family  $\{a_{I,J}:I,J \in \Ind_n\}$ which forms a magic unitary whose entries are selfadjoint projections, choose for each $I,J \in \Ind_n$ a subprojection $a_{I1,J1}$ of $a_{I,J}$ and define the projections $a_{I2,J1}, a_{I1,J2}, a_{I2,J2}$ by formulas \eqref{2final}, then the resulting family $\{a_{I',J'}:I',J' \in \Ind_{n+1}\}$ automatically yields a magic unitary whose entries are orthogonal projections.

Let us now see what the above discussion tells us about the structure of $\qg_n$ in this particular case. Start with $\alg_1$. As $\alg_0=\bc$, the invariance condition on the embedding  means simply that the action is unital and we see that $\qg_1$ is simply the universal quantum group acting on $2$ points, $C(\bz_2)$. Denote $a_{11} = p$. Then the unitary matrix coresponding to $\{a_{I,J}:I,J \in \Ind_1$ looks as follows:
\[ \left(\begin{array}{cc} p & p^{\perp} \\ p^{\perp}& p \end{array}\right).\]
In the second step, according to formulas \eqref{2final}, the 4 by 4 magic unitary looks as follows:
\begin{equation} \label{magun} \left(\begin{array}{cccc}  q_1 & p - q_1 & q_2 & p^{\perp} - q_2 \\
                                                            p - q_1 & q_1 & p^{\perp} - q_2 & q_2 \\
                                                            q_3 & p^{\perp} - q_3 & q_4 & p - q_4\\
                                                            p^{\perp} - q_3 & q_3 & p - q_4 & q_4
 \end{array}\right).\end{equation}
where $q_1,q_3$ are orthogonal subprojections of $p$ and $q_2,q_4$ are orthogonal subprojections of
$p^{\perp}$. The quantum group $\qg_2$ is then generated by the projections $p,q_1,q_2,q_3,q_4$
subjected to constraints described in the previous sentence (note that the subordination for
projections can be formulated in terms of the usual algebraic relations: $q\leq q'$ if and only if
$qq'=0$). $C^*$-algebraic structure of $\qg_2$ can be seen as follows - the projection $p$ provides
a decomposition of $\qg_2$ into a direct sum, and then in each of the factors we choose
independently two orthogonal projections (respectively $q_1$ and $q_3$ or $q_2$ and $q_4$). As
there are no relations between $q_1$ and $q_3$, the universal $C^*$-algebra they generate is simply
$C(\bz_2)\star C(\bz_2)$, the universal algebra generated by two orthogonal projections
(see \cite{Parthonote}). Thus
\[\qg_2 = \left(C(\bz_2)\star C(\bz_2)\right) \oplus  \left(C(\bz_2)\star C(\bz_2)\right).\]
 The Hopf$^*$- algebraic structure can be read immediately from the condition that the matrix \eqref{magun} gives a corepresentation of $\qg_2$.
       The inductive reasoning should now be clear. It can be visualised by the sequence of pictures, representing consecutive subdivisions of a square. The fractal structure of the limiting algebra is apparent. Note also that the classical symmetry group of the tree-type graph we consider can be graphically interpreted as a one-dimensional version of the above two-dimensional picture (so that the classical symmetry group $\qg^{\textrm{clas}}_n$ is simply equal to $\prod_{i=1}^{2^n} \bz_2$).
\end{proof}

We can also give a description of the $\lim_{n \in \bn} \qg_n$  in terms of the generators and relations:

\begin{cor}
The quantum isometry group of a natural spectral triple on the algebra of continuous functions on the middle-third Cantor set constructed in \cite{Chrivan} is the universal $C^*$-algebra generated by the family of selfadjoint projections
\[\{p\} \cup \bigcup_{n \in \bn} \left\{p_{m_1,\ldots m_n}:m_1,\ldots m_n \in \{1,2,3,4\}\right\}\]
subjected to the following relations:
\[p_1,p_2 \leq p, \;\;\; p_3, p_4 \leq p^{\perp},\]
\[ p_{m_1,\ldots,m_n, 1},  p_{m_1,\ldots,m_n, 2} \leq  p_{m_1,\ldots,m_n}, \;\;\;  p_{m_1,\ldots,m_n, 3},    p_{m_1,\ldots,m_n, 4}  \leq  p^{\perp}_{m_1,\ldots,m_n} \]
($n \in \bn, m_1,\ldots m_n \in \{1,2,3,4\}$).
\end{cor}
\begin{proof}
A straightforward consequence of Theorem \ref{induc} and Lemma \ref{mid3}.
\end{proof}

It is clear from the proof of the Theorem \ref{mid3} and also from the discussion before Lemma \ref{embed} that the quantum group actions we consider are actions on the tensor product of algebras preserving in some sense one of the factors. In the classical world such actions have to have a product form, as the following lemma confirms:

\begin{lem}
Suppose that $(X,d_X), (Y,d_Y)$ are compact metric spaces and $T:X\times Y \to X\times Y$ is an isometry satisfying the following condition: $\alpha_T (C(X) \ot 1_{Y}) \subset C(X) \ot 1_{Y}$, where $\alpha_T:C(X \times Y) \to C(X \times Y)$ is given simply by the composition with $T$. Then $T$ has to be a product isometry.
\end{lem}
\begin{proof}
Denote the family of isometries of $X\times Y$ satisfying the conditions of the lemma by $\isox$.
We claim that $\isox$ is a group. Recall that $\isoz$, the family of all isometries of a compact metric space $(Z,d_Z)$, is a compact group when considered with the topology of uniform convergence (equivalently, pointwise convergence; equivalently, metric topology given by $d (T_1,T_2) = \sum_{i=1}^{\infty} \frac{1}{2^i}d_Z(T_1(z_i), T_2(z_i))$, where $\{z_i:i \in \bn\}$ is a countable dense subset of $Z$). It is easy to that $\isox$ is a unital closed subsemigroup of $\isoyx$. Thus it is a compact semigroup satisfying the cancellation properties and it has to be closed under taking inverses.

Suppose now that $T\in \isox$. Then if $f\in C(X)$ we have for all $x\in X, y,y' \in Y$
\[ (f \ot 1_Y) (T(x,y)) = \alpha_T(f \ot 1_Y) (x,y) = \alpha_T(f \ot 1_Y) (x,y')  = (f \ot 1_Y) (T(x,y')).\]
This is equivalent to the fact that $T$  is given by the formula
\[T(x,y) = (h(x), g(x,y)), \;\;\; x\in X, y \in Y,\]
for some transformations $h:X \to X$, $g:X\times Y \to Y$. The fact that $T$ is an isometry implies in particular that for all $x,x' \in X$, $y \in Y$.
\begin{equation} \label{distance}d_X(x,x') = d_X(h(x), h(x')) + d_Y(g(x,y),g(x',y)).\end{equation}
In particular $h:X \to X$ is a contractive transformation. As by the first part of the proof
$T^{-1} \in \isox$, there exist transformations $h':X \to X$, $g':X\times Y \to Y$ such that
\[T^{-1}(x,y) = (h'(x), g'(x,y)), \;\;\; x\in X, y \in Y,\]
It is easy to see that $h'$ is the inverse transformation of $h$, and as by the same argument as
above we see that $h'$ is a contractive transformation, hence $h$ has to be an isometry. This together
with formula \eqref{distance} implies that $g:X \times Y \to Y$ does not depend on the first
coordinate, so that $T$ must be a product isometry. In particular $\isox = \mathsf{ISO}_X \times
\mathsf{ISO}_Y$.
\end{proof}

Theorem \ref{mid3} shows that the result above has no counterpart for quantum group actions, even on classical spaces.  We could have thought of elements of $\qg_2$ as quantum isometries acting on the Cartesian product of 2 two-point set, `preserving' the first coordinate in the sense analogous to the one in the lemma above. If this forced elements of $\qg_2$ to be product isometries, we would necessarily have $\qg_2 = \qg_1 \ot \qg_1$; in particular $\qg_2$ would have to be commutative.

\section{Quantum symmetry groups of finite metric spaces and finite graphs as quantum isometry groups of certain natural spectral triples}

Theorem \ref{AFgraph} shows that the quantum symmetry group of a particular type of a finite graph
(as defined in \cite{graph}) coincides with the quantum isometry group of a certain spectral
triple. Motivated by this we show in this section that given a finite metric space $X$ the quantum
symmetry group of $X$ defined in \cite{metric} coincides with the quantum isometry group of the
algebra of functions on $X$ equipped with the natural Laplacian. We also discuss the connections
with the natural Dirac operator on a particular representation of $C(X)$, for which the
Rieffel-type metric (\cite{Riefstate}) on the state space of $C(X)$ restricts to the original
metric on $X$. It can be seen, as pointed earlier in \cite{metric}, that this framework can be
related to the one of quantum symmetry groups of finite graphs.

Let $(X,d)$ be a finite metric space of $n$ points. For simplicity we will write $X=\{ 1,...,n \}$
and for $i,j  \in X$ define $d_{ij}=d(i,j)$. As in the  Section 2 we will denote by $\delta_i$
the indicator function of the point $i\in X$ and  by $ D_{ X \times X }$ the diagonal in $X \times
X$.

The metric structure on $X$ allows the construction of a natural spectral triple on $C(X)$ (see
\cite{Riefstate}, \cite{Chrivan2}). Let $ Y =  X \times X - D_{ X \times X },$ define the Hilbert
space $\Hil = \oplus_{( x_{0}, x_{1} ),x_{0}\neq x_{1}  } \Hil_{( x_{0}, x_{1} )},$ where $ \Hil_{(
x_{0}, x_{1} )} = {\bc}^{2},$ and let the Dirac operator be given by $$  D = \oplus_{( x_{0}, x_{1}
),x_{0}\neq x_{1}  } d^{-1}( x_{0}, x_{1} ) \left ( \begin {array} {cccc}
   0 & i  \\ -i & 0 \end {array} \right ) .$$
When we view $ \Hil $  as $ {\bc}^{2} \otimes l^{2} ( Y ),$ then \[ D = \left ( \begin {array} {cccc}
   0 & i  \\ -i & 0 \end {array} \right ) \otimes M_{d^{-1}}, \] where $ M_{d^{-1}} $ denotes multiplication by the function $ d^{-1} $ on $ Y .$

Let $ s, t:Y \rightarrow X $ be given by the formulas $ s( x_{0}, x_{1} ) = x_{0}, t ( x_{0}, x_{1} ) = x_{1}.$ Then for each $f
\in C(X) $ there is $\, s^{*} ( f ) = ( f \otimes 1 ) \chi_{Y},  t^{*} ( f ) = ( 1 \otimes f ) \chi_{Y}$, where $\chi_Y$ denotes
the characteristic function of $Y$. Let for $f \in C(X)$
\[\pi_{1} ( f ) =  \left ( \begin {array} {cccc}
   1 & 0  \\ 0 & 0 \end {array} \right ) \otimes s^{*} ( f ), \,\,\;
        \pi_{2} ( f ) =  \left ( \begin {array} {cccc}
   0 & 0  \\ 0 & 1 \end {array} \right ) \otimes t^{*} ( f ),\]
define  $ \pi( f ) = \pi_{1}( f ) + \pi_{2}( f )$ and consider the resulting spectral triple
$(C(X),  \Hil, D)$. Denote by $ {\bf C_{X,d}}$ the category of compact quantum groups acting by
volume (corresponding to $R = I$) and orientation preserving isometries on the above spectral
triple. By Theorem 2.10 in \cite{Deb2}, the universal object in this category exists; we denote it
by $\QI_I(D)$.

      It is also easy to describe the `Laplacian' in the sense of \cite{Deb} for this spectral triple.

\begin{lem}
  \label{laplacian}
  The Laplacian ${\mathcal L}$ on $C(X)$ associated with the spectral triple constructed above via the prescription in \cite{Deb} is given by the formula
\begin{equation} {\mathcal L}( \delta_{i} ) = \frac{4}{2n - 1} \sum_{j\in X} c ( i,j ) \delta_{j}, \;\;\;i \ in\ X\label{Lap} \end{equation} where for $i,j \in X$ we have  $c( i, j ) = d^{-2} ( i, j )$  if $i \neq j $ and $c(i,i)=0$.
 \end{lem}

  \begin{proof}

Let $ \tau_{0} $ be the functional on $C(X) $ defined by $ \tau_{0} ( f ) = \sum_{i = 1}^{n} f ( i
).$ Denote by $\Hil_0$ the Hilbert space obtained by  completing  $ \pi (C(X)) $ with respect to
the norm coming from the functional $ Tr(\pi(f)).$ It is easy to see that $ Tr ( \pi ( f ) ) = ( 2n
- 1 ) \tau_{0} ( f ).$

From the definition of the inner product it follows that for all $f,g \in C(X) $
\[\left\langle {\mathcal L } ( f ),
g  \right\rangle = \sum_{i\in X} g ( i ) ( 2n - 1 )  \overline{{\mathcal L} ( f )( i )}.\]
 On the other hand,
 $ \left\langle {\mathcal L} ( f ), g  \right\rangle = -\left\langle  d_D^{*} d_D f, g \right\rangle  =-
 \left\langle d_Df, d_Dg \right\rangle$ (where $d_D(\cdot)=[D,\cdot]$  as in \cite{Deb}),
 which by a routine calculation can be shown to be equal to
 $ 4 \sum_{i \neq j} \frac{g(i)(\overline{ f} ( j ) -\overline{ f} ( i) }{d^{2} ( i, j )}.$
 Thus for each $i \in X$
\[ {\mathcal L }( f ) ( i ) = \frac{4}{2n - 1} \sum_{j \neq i} \frac{f( j) - f( i ) }{d^{2} ( i, j )} .\]
Comparison of the above with formula \eqref{Lap} ends the proof.
\end{proof}

It is now easy to verify that $\mathcal L$ is admissible in the sense of \cite{Deb}, so that the
corresponding quantum isometry group  $\qil(X)$ exists. Recall that $\qil(X)$ is the universal
object in the category ${\bf C}^{\mathcal L}_{\bf X,d}$, with the objects being pairs $(\qg,
\alpha)$, where $\qg$ is a compact quantum group and $\alpha$ is the action of $\qg$ on $C(X)$
satisfying $({\mathcal L} \ot {\rm id}) \circ \alpha=\alpha \circ {\mathcal L}$.   We want to
compare $\qil(X)$  with the universal quantum symmetry group of $X$ constructed by Banica in
\cite{metric}. To this end we first need to observe the alternative characterisation of the actions
of compact quantum groups on finite metric spaces considered in \cite{metric}.

\begin{lem}
   \label{banica}
Given an action $ \alpha $ of a compact quantum group $ \qg $ on a finite metric space $(X,d)$
(i.e.\ an action of $\qg$ on $C(X)$), the following are equivalent:
\begin{rlist}
\item $(\qg, \alpha)$  is a quantum  isometry in the sense of Banica (\cite{metric});
\item   $ {\alpha}^{(2)}( d) = d \otimes 1$;
\item $\alpha^{(2)}(c)=c \ot 1$, where $c\in C(X \times X)$ is as in Lemma \ref{laplacian};
\item $(\qg, \alpha)$ is an object in ${\bf C}^{\mathcal L}_{\bf X,d}$.
\end{rlist}
\end{lem}
   \begin{proof}
(i) $\Leftrightarrow$ (ii):\\
   Write  $ d = \sum_{i,j\in X} d_{ij} \delta_{i} \otimes \delta_{j},$ where $d_{ij}$ are defined as in the beginning of this section. Let us write the action  $\alpha$ as
\begin{equation}\label{alpexp} \alpha(\delta_i)=\sum_{j\in X} \ \delta_j \ot q_{ij},\end{equation} where $q_{ij} \in \qg$. Then it follows by using the relations of the quantum permutation group that the relation $ {\alpha}^{(2)} (d) = d \otimes 1 $ is equivalent to the following equation being satisfied for all $k,l \in X$:
 \begin{equation} \label{qiso}  d_{kl} 1 = \sum_{i,j\in X} d_{ij} q_{ki} q_{lj}. \end{equation}
Thus, to prove the lemma, it is enough to show equivalence of  \eqref{qiso} with  Banica's
definition of quantum isometry.

Begin by noting that Banica's definition  implies that  for all $k,l \in X$
\[ \sum_{i\in X} d_{il} q_{ki} = \sum_{i\in X} d_{ki} q_{il} .\]
From this, it follows that for all $k,m \in X$,
\[\sum_{l,i\in X} d_{il} q_{ki} q_{ml} =  \sum_{l,i\in X} d_{ki} q_{il} q_{ml}
 =  \sum_{l \in X} d_{km} q_{ml} = d_{km}1, \]
which is exactly (\ref{qiso}).

For the converse direction rewrite (\ref{qiso}) as
 \begin{align*}  \sum_{i\in X} d_{ki} q_{il}
& =  \sum_{i\in X} ( \sum_{j,m\in X} d_{jm} q_{kj} q_{im} ) q_{il}
=  \sum_{i,j,m\in X} d_{jm} q_{kj} q_{im}  q_{il}\\&=
  \sum_{i,j,l\in X} q_{kj} \delta_{m,l} d_{jm} q_{im}  q_{il}~({\rm where}~ \delta_{m,l}~{\rm denotes~the~ Kronecker~delta})\\ &= \sum_{i,j\in X} q_{kj} d_{jl} q_{il} = \sum_{j\in X} d_{jl} q_{kj} . \end{align*}
Thus Banica's condition is derived. \vspace{1mm}\\
The proof of (iii) $\Leftrightarrow$ (iv) is very similar to the above proof of equivalence of (i)
and (ii), hence omitted.

Finally, the equivalence of (ii) and (iii) follows from the relation between $c$ and $d$, i.e.
$$ c = \chi_{Y} d^{-2},~~~    d = (1- \chi_{Y} + c )^{- \frac{1}{2}} - 1+\chi_{Y} ,$$  together with the fact that $\alpha^{(2)}(\chi_Y)=\chi_Y \ot 1$.
   \end{proof}

 Denote the category of compact quantum groups acting on the above spectral triple on a finite metric space
 $ (X,d)$ with action $ \alpha $ satisfying $ \alpha^{(2)} (d) =  d \otimes 1 $ by $ {\bf C}^{\bf Ban}_{\bf X,d} .$
 We want to show that that the universal object in $ {\bf C}^{\bf Ban}_{\bf X,d}$, say $\qib$
 (which is shown to exist in \cite{metric})  is isomorphic to the quantum group $\qi_I(D)$. The proof of
 this fact is contained in the two following lemmas.
 Recall that we always assume that $(X,d)$ is a finite metric space.

\begin{lem}
   \label{banicagoswami}
Let $ ( \qg, \alpha) $ be a quantum isometry of  $(X,d)$  in the sense of Banica, i.e.\ an object
in $ {\bf C}^{\bf Ban}_{\bf X,d}$. Then there is  a unitary representation $U$ of $\qg$ on $\Hil$
such that  $ (\qg, U) \in {\rm Obj}  ( { \bf C_{X,d} } ),$ with $\alpha_U=\alpha$ on $C(X)$.
      \end{lem}

\begin{proof}
      Define $ \widetilde{U} = I_{\bc^{2}} \otimes \widetilde{{\alpha}^{(2)}} $ on $ {\mathcal B} ( \Hil ) \otimes \qg.$ Then $\widetilde{U} $ gives a unitary representation since $ \widetilde{{\alpha}^{(2)}} $ is one.

 Moreover, recalling  from Lemma \ref{laplacian} that for all $f \in C(X)$ one has  $ Tr ( \pi ( f ) ) =  ( 2n - 1 ) \tau_{0} ( f )$ and that $ \alpha $ preserves $ \tau_{0},$ we immediately observe  that the action $ \alpha_U $ preserves the volume form corresponding to $ R = I$.
  Note that for any $X \in M_2(\bc)$
  \[ \widetilde{U} ( X \otimes M_{\phi} \otimes I_{\qg} ){ \widetilde{U}}^{-1} = X \otimes \widetilde{{\alpha}^{(2)}} ( M_{\phi} \otimes I ){ \widetilde{{\alpha}^{(2)}} }^{-1} = X \otimes \widetilde{ M_{\alpha^{(2)}( \phi )} },\] where for $Y=y_1 \otimes y_2 \in C(Y) \otimes \qg$, $\widetilde{M}_Y$ denotes $M_{y_1} \otimes y_2$.
Thus, in particular, \[ \widetilde{U} ( D \otimes I ){ \widetilde{U} }^{-1} = \left ( \begin {array} {cccc}
   0 & i  \\ -i & 0 \end {array} \right ) \otimes \widetilde{M_{{\alpha}^{(2)} (d^{-1})}} = \left ( \begin {array} {cccc}
   0 & i  \\ -i & 0 \end {array} \right ) \otimes M_{d^{-1}} \otimes I = D \otimes I,\]
where we have used $\alpha{(2)}(d^{-1})=d^{-1} \otimes 1$ which follows from  Lemma \ref{banica}. This  shows that $ U $ commutes with $ D .$
Further it is easy to see that  for any $f \in C(X)$,
\begin{align*}  \widetilde{U} ( \pi( f ) \otimes 1 ) {\widetilde{U}}^{-1}&=\\
 &= \widetilde{U} ( \pi_{1}( f ) \otimes 1 ) {\widetilde{U}}^{-1} + \widetilde{U} ( \pi_{2}( f ) \otimes 1 ) {\widetilde{U}}^{-1}\\
&=  \widetilde{U} \left( \left ( \begin {array} {cccc}
   1 & 0  \\ 0 & 0 \end {array} \right ) \otimes s^{*}( f )\otimes 1 \right) {\widetilde{U}}^{-1}
   + \widetilde{U} \left(  \left ( \begin {array} {cccc}
   0 & 0  \\ 0 & 1 \end {array} \right ) \otimes t^{*} ( f ) \otimes 1 \right){\widetilde{U}}^{-1} \\
&=   \left ( \begin {array} {cccc}
   1 & 0  \\ 0 & 0 \end {array} \right ) \otimes \widetilde{M_{{\alpha}^{(2)}( s^{*}f )}} +   \left ( \begin {array} {cccc}
   0 & 0  \\ 0 & 1 \end {array} \right ) \otimes \widetilde{M_{{\alpha}^{(2)}( t^{*}f )}} \\
&=   \left ( \begin {array} {cccc}
   1 & 0  \\ 0 & 0 \end {array} \right ) \otimes s^{*}( f_{(1)} ) \otimes f_{(2)} +   \left ( \begin {array} {cccc}
   0 & 0  \\ 0 & 1 \end {array} \right ) \otimes t^{*}( f_{(1)} ) \otimes f_{(2)} \\
&=  ( \pi \otimes \textup{id} ) \alpha ( f ) \subseteq \pi( C( X ) ) \otimes \qg. \end{align*} Note that in the above we have used the Sweedler notation $\alpha^{(2)}(f)=f_{(1)} \otimes f_{(2)}$. This
implies that  $ (\qg, U ) \in {\rm Obj} ( {\bf C_{X,d} }).$ It is obvious from the construction
that $\alpha_U=\alpha$.
   \end{proof}

 \begin{lem}
  \label{goswamibanica}
  Let $ (\widetilde{\qg}, U) \in {\rm Obj} ({\bf C_{X,d} }),$ with $\qg$ being
  the largest Woronowicz $C^*$-subalgebra of $\widetilde{\qg}$ such that the action
  $ \alpha_{U} $ maps $C(X)$ into $C(X) \ot \qg. $  Then  $ (\qg, \alpha_{U}) $ is an object of ${\bf C}^{\bf Ban}_{\bf X,d}$.
  \end{lem}

  \begin{proof}
 The fact that $ U $ commutes with $ D $ implies that $U$ commutes with $D^2=I_{\bc^2} \otimes M_{c}$ on $\bc^2 \otimes l^2(Y)$. Since $U=I \otimes \alpha^{(2)}$, it follows that $\alpha^{(2)}(c)=c \otimes 1$, hence (by Lemma \ref{banica})  $(\qg, \alpha)$ is a quantum isometry in the sense of Banica.
 \end{proof}

Lemma \ref{banica}, Lemma \ref{banicagoswami}  and Lemma \ref{goswamibanica} put together imply immediately the following:

 \begin{tw} \label{isom}
 We have the following isomorphisms of compact quantum groups: $$ \qib(X)  \cong  \qi_I(D) \cong \qil(X).$$
 \end{tw}



\begin{rem}
We can accommodate graphs in the framework of the above theorem if  we view a finite non-directed graph $(V,E)$ as a metric space $(V,d_E)$ where
\[d_E(v,w) = 1  \textrm{ if } (v,w) \in E, \;\;\;\; d_E(v,w) =  \infty \textrm{ if }  (v,w) \notin E \]
($v,w \in V, v \neq w$). A similar observation was made already in \cite{metric}. Here the Theorem \ref{isom} shows that quantum symmetry groups of finite graphs of \cite{graph} can be viewed as quantum isometry groups associated to the natural Laplacians on such graphs.
\end{rem}

\section{Remarks on quantum isometries of arbitrary metric spaces}

 In view of Lemma \ref{banica} it is reasonable to define a `quantum isometry' of a general metric space $(X,d)$
 to be a (faithful)  action $\alpha$ on $C(X)$ by a compact quantum group $\qg$ such that $\alpha^{(2)}(d)=d \ot 1$ (here again the metric $d$ is viewed as an element of $C(X \times X)$).  It is however a nontrivial (and open) problem to see whether there exists any universal objects in the category of such quantum isometries of $(X,d)$ when $X$ is not a finite set. We shall take up this issue elsewhere, but would like to
conclude this article with some computations for two simple yet interesting examples, namely the unit interval $[0,1]$ and the
circle $S^1$ (both equipped with the usual Euclidean metric).

\vspace*{0.2 cm}
$ \bf{ X = [ 0, 1 ] }$\vspace{1mm}\\
The $C^*$-algebra $C(X)$ is the universal unital $C^*$-algebra generated by a self-adjoint operator $T$ satisfying $0 \leq T \leq 1$. The metric is given by ($x_1, x_2 \in [ 0, 1 ]$)
\[ d^{2}( x_1 , x_2 ) = {x_{1}}^{2} - 2 x_{1}x_{2}  + {x_{2}}^{2}.\]
Thus, as an element of $ C ( [ 0, 1 ] ) \otimes C ( [ 0, 1 ] ),$ the metric $d$ is given by $ d^{2}
= T^2 \otimes 1 - 2 T \otimes T  + 1 \otimes T^2 $. Given a quantum isometric action $\alpha$ of a
c.q.g.\ $\qg$, let us write  \[ \alpha ( T ) =  \sum_{n \geq 0} T^n \otimes q_n \] (the series is
strongly convergent). Since  $ T $ is self adjoint and $ \alpha $ is a $^*$-homomorphism, each $
q_n $ is also self adjoint. Comparing coefficients of $ T^n \otimes T^n $ (for $n \geq 2$) and $
T^2 \otimes 1 $  in the equation $ \alpha^{(2)} (d^2) = d^2 \otimes 1 ,$ we get (respectively) $$
{q_{n}}^{2} =  0, \;\;\; n \geq 2 $$ and $$ q_{2}q_{0} + q_{0}q_{2} + q_{1}q_{1} - 2 q_{2}q_{0} =
1.$$
 Since  each $ q_{n} $ is self adjoint, the first equation above implies immediately that \[ q_{n} = 0, \;\;\; n \geq 2.\]
   Using $ q_2 = 0 $ in the second equation we deduce that  $ {q_{1}}^{2} = 1.$ Thus, $ q_1 $ is a reflection (selfadjoint unitary) and as such can be written as   $q_1 = P - P^{\bot},$ where $P$ is an orthogonal projection.

Moreover, $ \alpha ( T ) $ is a positive contraction as $ T $ is so, which implies that \[ 0 \leq 1
\otimes q_0 + T \otimes q_1 \leq 1.\] By applying to both sides of the above inequality the
evaluation functional at a point $t \in [0,1]$ we obtain
\begin{equation}\label{tineq} 0 \leq  q_0 + t q_1 \leq 1, \;\;\; t \in [ 0, 1 ].\end{equation}
Putting $ t = 0$ we note that $ q_0 $  is a positive contraction. Further
multiplying by $ P $ the both sides of \eqref{tineq}  we obtain \[ 0 \leq P q_0 P + t P \leq P, \;\;\;t \in [ 0, 1 ].\] Putting $ t = 1$ yields $ P q_0 P = 0 $ which implies that $ q_0 $ maps $ P $ to $ P^{\bot}$ (here and below we identify $P$ and $P^{\perp}$ with the corresponding subspaces of the Hilbert space).
Similarly, multiplying both sides of \eqref{tineq} by $ P^{\bot} $ we obtain  \[ 0 \leq P^{\bot} q_0 P^{\bot} - t P^{\bot} \leq P^{\bot}, \;\;\; t \in [0,1] .\]
Putting $ t = 0 $ and $t= 1,$ we get  (respectively)  $$ P^{\bot} q_0 P^{\bot} \leq P^{\bot} $$ and $$ P^{\bot} q_0 P^{\bot} \geq P^{\bot}.$$ It follows that $ P^{\bot} q_{0} P^{\bot} = P^{\bot} $, hence  $ q_0 : P^{\bot} \rightarrow P^{\bot} $ and $ q_0|_{P^{\bot}} = I .$
Recalling that $ q_0 $ is a contraction and using the above two observations we conclude that $ q_0 = P^{\bot} .$

Thus,  $ \alpha ( T ) = 1 \otimes P^{\bot} + T \otimes ( P - P^{\bot} ) $, which clearly implies  (by faithfulness of $\alpha$) that $\qg=C^*(P)$ is commutative, i.e. there is no `quantum isometry' of $X$. \vspace{1mm}\\

\vspace*{0.2 cm}
$ \bf{ X = S^{1} }$\vspace{1mm}\\
 The $C^*$-algebra $C(S^1)$ is the universal unital $C^*$-algebra generated by a unitary, say $Z$.
The metric is given this time by ($z_1, z_2 \in S^1$)
\[d^{2} ( z_1, z_2 ) = 2 - z_1 \overline{z_2} - z_2 \overline{z_1},\] so in the tensor picture we get
 $ d^2 = 2 \otimes 1 - Z \otimes Z^* - Z^* \otimes Z.$

Let $ \qg$ be a c.q.g.\ with an isometric action $ \alpha $ on $S^1$. Let \[\alpha ( Z ) = \sum_{n
= 0}^{\infty} Z^n \otimes q_n + \sum_{n = 1}^{\infty} {Z^{*}}^{n} \otimes q^{\prime}_{n} ,\] where
$ q_n, q^{\prime}_{n} \in \qg.$ Comparing the coefficients of  $ 1 \otimes 1 $  and $ Z^n \otimes
{Z^{*}}^{n}$  for $n \geq 2 $ on both sides of the equation $ \alpha^{(2)} (d^2) = d^2 \otimes 1 $
we get this time $$ q_0 q^{*}_0 + q^{*}_{0} q_{0} = 0, \;\;\;q_n q^{*}_n +  {q^{\prime }}^*_n
q^{\prime}_n = 0, \;\;\; n \geq 2. $$ This proves that  $ q_0 = 0 $ and $ q_n = q^{\prime}_n = 0$
for all $n \geq 2 $. Thus, $ \alpha ( Z ) = Z \otimes q_1 + Z^* \otimes q^{\prime}_{1} ,$
 i.e.  $ \alpha $ is  `linear' as in subsection 2.2 in  \cite{DebJyot} and hence by Theorem 2.4 of that paper,
 $\qg$ must be commutative as a $C^*$-algebra and so a quantum subgroup of $ C(ISO(X))=  C( S^1 >\!\!\! \lhd Z_{2} )  .$
Then one can conclude that the investigated quantum isometry group $ QISO(S^1, d)$  is equal to $C( S^1 >\!\!\! \lhd Z_{2} ) $,
which can be interpreted as a statement that there is no `quantum isometry' of $S^1$ viewed as a metric space. We remark that we
arrived at a similar conclusion in \cite{DebJyot} and \cite{Deb} by viewing $S^1$ as a Riemannian manifold. So, the observation
made in the present paper in some sense strengthens the results of \cite{DebJyot} and \cite{Deb} about the quantum actions on
$S^1$.

\vspace*{0.2cm} \noindent \textbf{Acknowledgment.} The main results of this paper have been obtained during the visit of the
third named author to the Indian Statistical Institute in Kolkata, which was made possible thanks to the support of the UKIERI
project Quantum Probability, Noncommutative Geometry and Quantum Information.

The first and second named authors would like to acknowledge support from the National Board of Higher Mathematics, D.A.E. (India) and Indian National Science Academy respectively.

\end{document}